\documentclass[11pt]{amsart}

\usepackage{apacite}

\usepackage{amsmath, amscd, amsthm, amssymb, mathrsfs,amsfonts}
\usepackage{graphicx}

\theoremstyle{definition}

\theoremstyle{remark}

\numberwithin{equation}{section}

\begin{document}

\title{The antinomy of M. Gödel}


\author{Translation by Jérôme Verstrynge}

\begin{abstract}
	In order to facilitate historical research in logic, we provide an english translation of a french article by Chaïm Perelman, called 'L'antinomie de M. Gödel', published in the Bulletins de l'Académie royale de Belgique in 1936. 
\end{abstract}

\maketitle




\begin{center}
\textbf{The antinomy of M. Gödel,}
\end{center}
\vspace{2px}
\begin{center}
by Ch. Perelman, F.N.R.S. aspirant
\end{center}

\vspace{10px}

In a now famous article\footnote{\phantom{ }Ueber formal unentscheidbare Sâtze der Principia Mathematica und verwandter Systeme I. Aus den \textit{Monatsheften für Mathematik und Physik}, XXXVIIIe, v., 1er cahier, Leipzig, 1931.}, M. Gödel has proposed to demonstrate that inside an analytical system, like that of \textit{Principia Mathematica}\footnote{\phantom{ }A. WHITEHEAD and B. RUSSELL. Principia Mathematica, 2nd ed., Cambridge, 1925.}, it is possible to construct indeterminable propositions.

By \textit{indeterminable} proposition, we mean a proposition such that one cannot demonstrate either the truth of the affirmation, or the truth of the negation of this proposition. In order to demonstrate that a proposition $p$ is indeterminable, one must reduce the hypothesis that $p$ is provable and that not-$p$ is provable to a contradiction. To settle some terminology, let's say that a \textit{provable} expression is true; by transposition, a false proposition is \textit{unprovable}.
	
In order to construct an indeterminable proposition inside an analytical system, M. Gödel arithmetizes the formulas of logic, that is to say, he makes a determined natural number correspond to each logical sign; so that each formula can be expressed using a finite sequence of natural numbers, and each demonstration using a finite sequence of finite sequences of natural numbers. Proving a proposition amounts to constructing a sequence of sequences of natural numbers such that the last sequence expresses the demonstrated proposition. Mr. Gödel shows that it is possible to construct, with the help of \textit{Principia Mathematica} symbols, a formula $F(v)$ whose signification would be "$v$ is a provable formula".

This being admitted, the construction of an indeterminable proposition turns out to be relatively easy.

Consider a sign corresponding to a propositional function; it will consist of a sequence of natural numbers containing a variable, that is, an empty place, where one can put a random argument, which will be a natural number. Since all the formulas of this calculation consist of finite sequences of numbers, the number of unique formulas one could construct cannot exceed the order of what is countable, and it is possible to number these formulas, that is to associate each of them with an index which will be a natural number.

If, amongst all the functions we have numbered, one considers as a value for the variable the numerical sign corresponding to the index of each of them, one will obtain sequences of natural numbers where some will be provable and others will not. This fact makes it possible to distinguish two kinds of indices: those linked to functions giving rise to provable expressions, when one gives the variable their index as a value, and those which do not have this property. \textit{The latter indexes can be grouped into a set $E$}. Saying "$n$ is an element of $E$" will therefore be equivalent to "it is false that by giving the variable of the function $_n Fx$ the value $n$ we obtain a provable expression".

Note that the propositional function "$n$ is element of E" also has an index, say $q$; one can therefore designate it with symbol $_q Fx$. By introducing this symbol into the last equivalence, we have "$n$ satisfies $_q Fx$" is equivalent to "it is false that $n$ \textit{satisfies} $_n Fx$ is a provable expression”, or again: "$_q Fn$ is equivalent to "$_n Fn$ is an unprovable expression".

From this last equivalence, it is easy to demonstrate that $_q Fq$ is an indeterminable expression, that is, the hypothesis of its provability, as much as the provability of its negation, results in a contradiction.

Indeed, suppose that $_q Fq$ is provable; in this case, $_q Fq$ will be a true expression, and it will be the same of its equivalent expression: "$_q Fq$ is an unprovable expression". Hypothesizing the provability of $_q Fq$ implies thus its negation, and must as a consequence be rejected.

Suppose that not-$_q Fq$ is provable; it will therefore be true, like its equivalent expression "it is false that $_q Fq$ is an unprovable expression"; but this is the same as claiming the demonstrability of $_q Fq$. Since hypothesizing the provability of not-$_q Fq$ also implies its negation, the $_q Fq$ expression must be considered as indeterminable.

Let's take a closer look at the reasoning of M. Gödel. A formal study of the assumptions supporting it will convince us that the result of M. Gödel is more modest than it seems at first glance. We will see, indeed, that the only result of the reviewed article is the construction of a new antinomy, to be added to those, now classic, of logic and of set theory. The antinomy of M. Gödel presents exactly the same structure as the latter, and results, like them, from a false equivalence laid down in premises\footnote{\phantom{ }V. CH. PERELMAN, Les paradoxes de la logique, in \textit{Mind }, april 1936, vol. XLV. For an english translation, see: https://hal.science/hal-04294226}.

All of Mr. Gödel's reasoning revolves around a fundamental definition, that of set E. An index is part of it if the expression, obtained by replacing the variable of the function by the index of the latter, is unprovable. This allows us to lay down the equivalence defining set E:

\begin{center}
	\vspace{3px}
	$(n) \cdot n \, \epsilon \, E \, \cdot \equiv \cdot \sim  Dem \, _n Fn.$
	\vspace{3px}
\end{center}

($\sim \cdot \, Dem \, _n Fn$ is to be read: it is false that $_n Fn$ is demonstrable.)

Since the $n \, \epsilon \, E$ expression designates a function whose index is $q$, we obtain, with a simple replacement in the previous equivalence,

\begin{center}
		\vspace{3px}
	$(n) \cdot \, _q Fn \, \cdot \equiv \cdot \sim \cdot Dem \, _n Fn.$
		\vspace{3px}
\end{center}

This equivalence being affirmed for all values of $n$, it must also be true if one replaces $n$ by one of its value, namely $q$:

\begin{center}
		\vspace{3px}
	$_q Fq \, \cdot \equiv \cdot \sim \cdot \, Dem \, _q Fq$;
		\vspace{3px}
\end{center}

however, this last formula affirms the equivalence of an expression with the affirmation of its unprovability. One can notice the analogy of this result with the liar paradox, where one affirms, in an expression, the falsity of it. Nonetheless, Mr. Gödel seems to bring off, not a paradox, but a valuable result. Indeed, while in paradoxes one arrives at the equivalence of a proposition and its negation, one has managed here to demonstrate the equivalence of a proposition, not with its negation, but with the affirmation of its unprovability. Now, this result does not seem, a priori, contradictory; instead of classifying it among the paradoxes, it was considered as a mathematical discovery of the highest importance.

In the following pages we will demonstrate that from the fundamental equivalence from which Mr. Gödel starts it is possible to deduce paradoxes laying down the equivalence of a proposition and its negation, that Mr. Gödel's reasoning is based on false premises, just like the other antinomies, and that one should therefore not be surprised to see Mr. Gödel achieve an extraordinary result, whereas it is easy to come, from these same premises, to a contradictory result.

Consider the equivalence:

\begin{center}
		\vspace{3px}
	$(n) \cdot n \, \epsilon \, E \, \cdot \equiv \cdot \sim \cdot \, Dem \, _n Fn.$ \phantom{  } (1)
		\vspace{3px}
\end{center}

Let's apply the transposition principle to it:

\begin{center}
		\vspace{3px}
	$(n) \, \cdot \sim n \, \epsilon \, E \, \cdot \equiv \cdot \, Dem \, _n Fn.$ \phantom{  } (2)
		\vspace{3px}
\end{center}

In what follows, we will transform these first equivalences so as to make their second member equivalent, not to the
truth of the first, but to its provability. We will start by transforming equivalence (2).

Let's note that set E has been defined as the one of all $n$, such that $\sim \cdot \, Dem \, _n Fn$. Its complementary will be the set of all $n$, that $Dem \, _n Fn$. It follows that if, for a given $n$, $_n Fn$, is provable, the expression $\sim n \, \epsilon \, E$ is not only true, but also provable; it will be enough, in fact, to complete the demonstration of $_n Fn$ by the definition of set E to demonstrate $\sim n \, \epsilon \, E$.

The way in which the set E was introduced allows us to lay down:

\begin{center}
		\vspace{3px}
	$(n) \cdot Dem \, _n Fn \;\supset\; Dem \, \sim  n \, \epsilon \, E.$ \phantom{  } (3)
		\vspace{3px}
\end{center}

On the other hand, since every provable expression is true, we have

\begin{center}
		\vspace{3px}
	$(n) \cdot Dem \, \sim n \, \epsilon \, E \;\supset\; \sim n \, \epsilon \, E.$ \phantom{  } (4)
		\vspace{3px}
\end{center}

From equivalence (2), we can get:

\begin{center}
		\vspace{3px}
	$(n) \, \cdot \sim n \, \epsilon \, E \;\supset\; Dem \, _n Fn.$ \phantom{  } (5)
		\vspace{3px}
\end{center}

Since from (4) and (5), one can get, by syllogism,

\begin{center}
		\vspace{3px}
	$(n) \cdot Dem \, \sim n \, \epsilon \, E \;\supset\; Dem \; _n Fn$, (6)
		\vspace{3px}
\end{center}

one can deduce from (3) and (6)

\begin{center}
		\vspace{3px}
	$(n) \cdot Dem \, \sim n \, \epsilon \, E \, \cdot \, \equiv \, \cdot \, Dem \, _n Fn$. \phantom{  }  (7)
		\vspace{3px}
\end{center}

Since $n \, \epsilon \, E$ can be written $_q Fn$, proposition (7) is equivalent to

\begin{center}
		\vspace{3px}
	$(n) \cdot Dem \, \sim \, _q Fn \, \cdot \equiv \cdot \, Dem \, _n Fn$. \phantom{  } (8)
		\vspace{3px}
\end{center}

By applying the same reasoning to proposition (1) which resulted in proposition (8), using equivalence (2), we obtain the equivalence

\begin{center}
		\vspace{3px}
	$(n) \cdot Dem \, _q Fn \, \cdot \equiv \cdot \sim \cdot \, Dem \, _n Fn$. (9)
		\vspace{3px}
\end{center}

For argument $q$, equivalences (8) and (9) become

\begin{center}
		\vspace{3px}
	$Dem \, \sim \, _q Fq \, \cdot \equiv \cdot \, Dem \, _q Fq$ \phantom{  } (10)
		\vspace{1px}
\end{center}

and 

\begin{center}
		\vspace{1px}
	$Dem \, _q Fq \, \cdot \equiv \cdot \sim \cdot \, Dem \, _q Fq$. \phantom{  } (11)
		\vspace{3px}
\end{center}

Considering propositions (10) and (11), which one could easily deduce from the definition of set E, laid down by Mr. Gödel, one will not be surprised by the result obtained by him, which exploits only to a small extent the possibilities it offered. It makes it possible, in fact, to demonstrate, in addition to the existence of an indeterminable proposition, a proposition whose truth and falsity can be demonstrated simultaneously, and a proposition whose provability is equivalent to unprovability. Here we are, right in the middle of paradoxes. This harvest of results, however, is richer than that obtained in classical paradoxes where one is satisfied enough to demonstrate that the truth of a proposition is equivalent to its falsity; in fact, the notion of provability allows us to obtain four possibilities for a given proposition: it can be provable or unprovable, and the same is true of its negation.

How was it possible to demonstrate propositions (10) and (11), whose contradictory character is nevertheless explicit? They were deduced from propositions (8) and (9), which constitute false equivalences. The first states the equivalence of the provability of a function to that of its negation, and this for all the values of the variable; in the second proposition, we affirm the formal equivalence of the provability of a function to its unprovability. The fact that in the first member of these equivalences is located, in place of the variable, one of the arguments of the function makes it possible to limit the contradiction; but it explodes when we introduce the same argument into the second member of the equivalence.

As these false equivalences are deduced, with the greatest ease, from the definition of set E, it is this which must be called into question; by taking a close look, in fact, one realizes that it consists of a false equivalence, that the principle of contradiction forces us to reject.

By saying, in his article, that his demonstration was related to paradoxes, Mr. Gödel said too little. It is, in fact, a new antinomy, having identically the same structure as known paradoxes, which he had just build. And, just like classic paradoxes, the one that Mr. Gödel constructed in his famous article results from a contradiction lying in the premises.		
		

\end{document}